\documentclass[11pt,twoside]{article}
\newcommand{\ifims}[2]{#1}   
\newcommand{\ifAMS}[2]{#1}   
\newcommand{\ifau}[3]{#1}  

\newcommand{\ifbook}[2]{#1}   

\ifbook{

  }
  {

  }

\def\thetitle{A note on critical dimensions in profile semiparametric estimation}
\def\theruntitle{A note on  critical dimensions in semiparametric estimation}

\def\theabstract{
This paper complements the results of \citep{AASP2013} on profile estimators in semiparametric models. We present two examples. One that illustrates that the smoothness constraint on the expected value of the contrast functional used to define the profile M-estimator is necessary for the bound derived for the critical ratio of dimension to sample size. A second one to show that in the case that the target dimension is proportional to the full dimension the critical ratio for the Fisher type result stays the same while for the Wilks phenomenon it is multiplied with the square root of the full dimension, just as in the upper bound in \citep{AASP2013}.
} 

\def\kwdp{62F10}
\def\kwds{62J12,62F25,62H12}

\def\thekeywords{profile, semiparametric, spread,
local concentration, critical dimension}

\def\authora{Andreas Andresen }
\def\runauthora{andreas andresen }
\def\addressa{
    \\ Weierstrass-Institute, \\
    Mohrenstr. 39, \\
    10117 Berlin, Germany     
    }
\def\emaila{andresen@wias-berlin.de}
\def\affiliationa{}
\def\thanksa{The author is supported by Research Units 1735 
"Structural Inference in Statistics: Adaptation and Efficiency"
}

\def\authorb{Vladimir Spokoiny }
\def\runauthorb{vladimir spokoiny }
\def\addressb{
    Weierstrass Institute and HU Berlin, \\ Moscow Institute of 
    Physics and Technology \\
    Mohrenstr. 39, \\
    10117 Berlin, Germany}
\def\emailb{spokoiny@wias-berlin.de}
\def\affiliationb{Weierstrass-Institute, Humboldt University Berlin, and
Moscow Institute of Physics and Technology}
\def\thanksb
{The author is partially supported by
Laboratory for Structural Methods of Data Analysis in Predictive Modeling, MIPT, 
RF government grant, ag. 11.G34.31.0073.
Financial support by the German Research Foundation (DFG) 
through the CRC 649 ``Economic Risk'' is gratefully acknowledged. 
}

\date{}

\usepackage{amsmath,amssymb,amsthm}
\usepackage{natbib}
\usepackage{epsfig,graphicx}
\usepackage{comment}
\usepackage{color}
\usepackage{srcltx}
\usepackage[mathscr]{eucal}
\usepackage[math]{easyeqn}
\usepackage{etoolbox}

\ifims{
\textheight=23cm
\textwidth=14.8cm
\topmargin=0pt
\oddsidemargin=1.0cm
\evensidemargin=1.0cm
\linespread{1.3}
\renewenvironment{abstract}
    {\centerline{\textbf{Abstract}}\bigskip
      \begin{center}
       \begin{minipage}{11cm}
        \begin{small}
    }
    {   \end{small}
       \end{minipage}
      \end{center}
     \bigskip
    }

}{ 
}

\numberwithin{equation}{section}
\numberwithin{figure}{section}
\newcounter{example}[section]
\numberwithin{example}{section}
\newcounter{remark}[section]
\numberwithin{remark}{section}
\newtheorem{theorem}{Theorem}[section]

\newtheorem{lemma}[theorem]{Lemma}

\newtheorem{exmp}[example]{Example}
\newtheorem{rmrk}[remark]{Remark}
\newenvironment{example}{\begin{exmp}\rm}{\end{exmp}}
\newenvironment{remark}{\begin{rmrk}\rm}{\end{rmrk}}

\begin{document}
\thispagestyle{empty}
\ifims{
\title{\thetitle}
\ifau{ 
  \author{
    \authora
    \ifdef{\thanksa}{\thanks{\thanksa}}{}
    \\[5.pt]
    \addressa \\
    \texttt{ \emaila}
  }
}
{  
  \author{
    \authora
    \ifdef{\thanksa}{\thanks{\thanksa}}{}
    \\[5.pt]
    \addressa \\
    \texttt{ \emaila}
    \and
    \authorb
    \ifdef{\thanksb}{\thanks{\thanksb}}{}
    \\[5.pt]
    \addressb \\
    \texttt{ \emailb}
  }
}
{   
  \author{
    \authora
    \ifdef{\thanksa}{\thanks{\thanksa}}{}
    \\[5.pt]
    \addressa \\
    \texttt{ \emaila}
    \and
    \authorb
    \ifdef{\thanksb}{\thanks{\thanksb}}{}
    \\[5.pt]
    \addressb \\
    \texttt{ \emailb}
    \and
    \authorc
    \ifdef{\thanksc}{\thanks{\thanksc}}{}
    \\[5.pt]
    \addressc \\
    \texttt{ \emailc}
  }
}

\maketitle
\pagestyle{myheadings}
\markboth
 {\hfill \textsc{ \small \theruntitle} \hfill}
 {\hfill
 \textsc{ \small
 \ifau{\runauthora}
      {\runauthora and \runauthorb}
      {\runauthora, \runauthorb, and \runauthorc}
 }
 \hfill}
\begin{abstract}
\theabstract
\end{abstract}

\ifAMS
    {\par\noindent\emph{AMS 2000 Subject Classification:} Primary \kwdp. Secondary \kwds}
    {\par\noindent\emph{JEL codes}: \kwdp}

} 
{ 
\begin{frontmatter}
\title{\thetitle}


\runtitle{\theruntitle}

\ifau{ 
\begin{aug}
    \author{\authora\ead[label=e1]{\emaila}}
    \address{\addressa \\
     \printead{e1}}
\end{aug}

 \runauthor{\runauthora}
\affiliation{\affiliationa} }
{ 
\begin{aug}
    \author{\authora\ead[label=e1]{\emaila}\thanksref{t21}}
    \and
    \author{\authorb\ead[label=e2]{\emailb}\thanksref{t22}}
    
    \address{\addressa \\
     \printead{e1}}
    \address{\addressb \\
     \printead{e2}}
    \thankstext{t21}{\thanksa}
    \thankstext{t22}{\thanksb}
    \affiliation{\affiliationa, \affiliationb} 
    \runauthor{\runauthora and \runauthorb}
\end{aug}
} 
{ 
\begin{aug}
    \author{\authora\ead[label=e1]{\emaila}\thanksref{t21}}
    \and
    \author{\authorb\ead[label=e2]{\emailb}\thanksref{t22}}
    \and
    \author{\authorc\ead[label=e3]{\emailc}\thanksref{t23}}
    
    \address{\addressa \\
     \printead{e1}}
    \address{\addressb \\
     \printead{e2}}
    \address{\addressc \\
     \printead{e3}}
    \thankstext{t21}{\thanksa}
    \thankstext{t22}{\thanksb}
    \thankstext{t23}{\thanksc}
    \affiliation{\affiliationa, \affiliationb, \affiliationc} 
    \runauthor{\runauthora, \runauthorb, and \runauthorc}
\end{aug}}

\begin{abstract}
\theabstract
\end{abstract}



\end{frontmatter}
} 

\def\ND{\cc{N}}
\def\Bernoulli{\mathrm{Bernoulli}}
\def\Vola{\mathrm{Vola}}
\def\Poisson{\mathrm{Poisson}}
\def\ag{\mathrm{ag}}
\def\glob{\operatorname{glob}}
\def\blk{\operatorname{block}}
\def\lin{\operatorname{lin}}
\def\cond{\, \big| \,}

\def\rdl{\epsilon}
\def\rd{\bb{\rdl}}
\def\rddelta{\delta}
\def\rdomega{\varrho}
\def\rddeltab{\rddelta^{*}}
\def\rhorb{\rhor^{*}}

\def\wv{\bb{w}}
\def\varthetav{\bb{\vartheta}}
\def\Lr{\breve{L}}
\def\zetavr{\breve{\zetav}}
\def\etavr{\breve{\etav}}
\def\xivr{\breve{\xiv}}

\def\rdb{\rd}
\def\rdm{\underline{\rdb}}

\def\taub{\tau_{\rdb}}
\def\taum{\tau_{\rdm}}
\def\kappab{\kappa_{\rd}}
\def\deltab{\delta_{\rd}}

\def\taubGP{\tau_{\rdb,\GP}}
\def\taumGP{\tau_{\rdm,\GP}}
\def\kappabGP{\kappa_{\rd,\GP}}
\def\deltabGP{\delta_{\rd,\GP}}
\def\nubm{\nu_{\rd}}
\def\uub{u_{\rd}}
\def\uubGP{u_{\rd,\GP}}
\def\nubmGP{\nu_{\rd, G}}

\def\rG{\rd,\GP}

\def\LinSp{\mathrm{L}}
\def\Id{I\!\!\!I}
\def\Ind{\operatorname{1}\hspace{-4.3pt}\operatorname{I}}

\def\BG{\mathcal{R}}
\def\bg{r}
\def\fmup{\phi}
\def\rg{r}
\def\uc{u_{c}}
\def\muc{\mu_{c}}
\def\mud{\mu_{0}}
\def\xxd{\xx_{0}}
\def\yyd{\yy_{0}}
\def\gmd{\gm_{0}}

\def\ms{m^{*}}
\def\Inv{A}
\def\InvT{\Inv^{\T}}
\def\Invt{\tilde{\Inv}}

\def\ssize{N}
\def\nsize{{n}}

\def\rhor{\omega}

\def\LT{L}
\def\LGP{\LT_{\GP}}
\def\La{\mathbb{L}}
\def\Lab{\La_{\rdb}}
\def\Lam{\La_{\rdm}}

\def\DP{D}
\def\DPc{\DP_{0}}
\def\DPb{\DP_{\rdb}}
\def\DPm{\DP_{\rdm}}

\def\LabGP{\La_{\rdb,\GP}}
\def\LamGP{\La_{\rdm,\GP}}

\def\DPbGP{\DP_{\rdb,\GP}}
\def\DPmGP{\DP_{\rdm,\GP}}
\def\riskbGP{\riskt_{\rdb,\GP}}

\def\gmi{\mathtt{b}}
\def\gmiid{\mathtt{g}_{1}}
\def\kullbi{\Bbbk}
\def\Thetasi{\Theta_{\loc}}
\def\rri{\mathtt{u}}
\def\rris{\rri_{0}}

\def\Ipc{\bb{\mathrm{f}}}
\def\IF{\Bbb{F}}
\def\IFc{\IF_{0}}
\def\IFb{\IF_{\rdb}}
\def\IFm{\IF_{\rdm}}

\def\DF{\cc{D}}
\def\DFc{\DF_{0}}
\def\DFb{\DF_{\rdb}}
\def\DFm{\breve{\DF}_{\rd}}
\def\DFm{\DF_{\rdm}}

\def\DPr{\breve{\DP}}
\def\VF{\cc{V}}
\def\VFc{\VF_{0}}

\def\HHc{\HH_{0}}
\def\HHb{\HH_{\rd}}
\def\HHm{\HH_{\rdm}}

\def\xib{\xi^{*}}
\def\xivb{\xiv_{\rdb}}
\def\xivm{\xiv_{\rdm}}
\def\CAm{\underline{\CA}}
\def\CAb{\CA}

\def\penr{\operatorname{pen}}
\def\pen{\mathfrak{t}}
\def\PEN{\operatorname{PEN}}
\def\RSS{\operatorname{RSS}}
\def\med{\operatorname{med}}

\def\ex{\mathrm{e}}
\def\entrl{\mathbb{Q}}
\def\entrlb{\entrl}
\def\entr{\entrl}

\def\kullb{\cc{K}} 
\def\kullbc{\kullb^{c}}

\def\gm{\mathtt{g}}
\def\gmc{\gm_{c}}
\def\gmb{\gm}
\def\gmbm{\gmb_{1}}

\def\yy{\mathtt{y}}
\def\yyc{\yy_{c}}
\def\xx{\mathtt{x}}
\def\xxc{\xx_{c}}
\def\tc{t_{c}}

\def\alp{\alpha}
\def\alpn{\rho}
\def\gmu{\mathfrak{a}}

\def\losst{\varrho}
\def\loss{\wp}
\def\lossp{u}
\def\closs{g}

\def\riskt{\cc{R}}
\def\emprisk{\ell}
\def\bias{b}
\def\bern{q}

\def\TT{\nsize}

\def\Pone{P}
\def\Pf{\P_{f(\cdot)}}
\def\Ef{\E_{f(\cdot)}}
\def\Ps{\P_{\thetas}}
\def\Es{\E_{\thetas}}
\def\Pu{\P_{\upsilons}}
\def\Eu{\E_{\upsilons}}

\def\Pvs{\P_{\thetavs}}
\def\Evs{\E_{\thetavs}}

\def\UPd{w}
\def\nunup{\nu_{1}}
\def\rru{\rr_{1}}
\def\rups{\rr_{0}}
\def\rupsb{\rups^{*}}
\def\rrf{\rr^{\flat}}
\def\rupd{\rr_{\circ}}

\def\smooths{\mathbb{S}}
\def\smooth{\smooths_{1}}

\def\elli{\bar{\ell}}

\def\K{K}

\def\Psir{\breve{\Psi}}

\def\af{a}
\def\afs{\af^{*}}

\def\kapla{\varkappa}

\newcommand{\mlew}[1]{\tilde{\thetav}_{#1}}
\newcommand{\mlea}[1]{\hat{\thetav}_{#1}}
\newcommand{\mluw}[1]{\tilde{\theta}_{#1}}
\newcommand{\mlua}[1]{\hat{\theta}_{#1}}
\newcommand{\penm}[1]{\boldsymbol{m}_{#1}}

\def\Pdom{\mu_{0}}
\def\PDOM{\bb{\mu}_{0}}
\def\EDOM{\E_{0}}

\def\mk{m}
\def\Mk{\cc{M}}
\def\SV{\cc{S}}

\def\Cs{E}
\def\Csd{\Cs^{\circ}}
\def\Ca{A}
\def\CS{\cc{E}}
\def\CA{\cc{A}}
\def\CAb{\CA_{\rd}}
\def\CAC{\CA_{\CoFu}}

\def\Ccb{m_{\rdb}}
\def\Ccm{m_{\rdm}}
\def\CcbGP{m_{\rdb,\GP}}
\def\CcmGP{m_{\rdm,\GP}}

\def\etas{\eta^{*}}

\def\omegav{\bb{\phi}}
\def\omegavs{\omegav^{*}}
\def\omegavc{\omegav'}

\def\nuvs{\nuv^{*}}
\def\nuvc{\nuv'}

\def\nunu{\nu_{0}}
\def\numu{\nu_{1}}
\def\nupi{\nu^{+}}
\def\nubu{\beta}

\def\nus{\nu}
\def\nusb{\nus}
\def\nusr{\nus^{\bracketing}}
\def\Nusb{\mathbb{N}}
\def\Nusr{\mathbb{N}^{\diamond}}

\def\dist{d}
\def\distd{\mathfrak{a}}

\def\hatk{\kappa}
\def\ko{k^{\circ}}

\def\qqq{\mathfrak{q}}
\def\ppp{{s}}
\def\Cqq{C(\qqq)}
\def\Cqqb{C^{\diamond}(\qqq)}
\def\Crho{C(\mrho)}
\def\Cqqm{\log(4)}
\def\Cqpr{(\qqq \rrp + \dimp / 2)}

\def\Cdima{\mathfrak{e}_{0}}
\def\Cdimb{\mathfrak{e}_{1}}
\def\cdima{\mathfrak{c}_{0}}
\def\cdimb{\mathfrak{c}_{1}}
\def\cdim{\mathfrak{c}}

\def\rdomega{\varrho}
\def\deltaD{\delta}
\def\alphai{\alpha_{1}}
\def\alphaii{\alpha_{2}}
\def\alphaiii{\alpha_{3}}
\def\alphaiv{\alpha_{4}}

\def\err{\diamondsuit}
\def\errbm{\bar{\err}_{\rdomega}}
\def\errm{\err_{\rdm}}
\def\errb{\err_{\rdb}}

\def\errbGP{\err_{\rdomega,\GP}}
\def\errmGP{\err_{\rdm,\GP}}
\def\errbmGP{\bar{\err}_{\rd,\GP}}

\def\errs{\err_{\rdomega}^{*}}
\def\deltas{\alpha}

\def\xivbGP{\xiv_{\rdb,\GP}}
\def\xivmGP{\xiv_{\rdm,\GP}}

\def\SP{S}
\def\GP{G}
\def\GPt{\GP_{0}}
\def\GPn{\GP_{1}}
\def\gp{g}
\def\gs{s}

\def\errbGP{\err_{\rdb,\GP}}
\def\errmGP{\err_{\rdm,\GP}}
\def\errpmGP{\err_{\GP}^{\pm}}

\def\LCS{\cc{C}}

\def\DPGP{\DP_{\GP}}
\def\thetavsGP{\thetavs_{\GP}}

\def\LL{\cc{L}}
\def\LLb{\LL^{*}}
\def\LLh{\cc{L}}

\def\YY{\cc{Y}}
\def\LP{L^{\circ}}

\def\modcnrd{\mathfrak{A}}

\def\pens{\pi}
\def\pnn{\mathfrak{g}}
\def\pnnd{\mathfrak{u}}
\def\pnndGP{\pnnd_{\GP}}

\def\confpr{\mathfrak{c}}
\def\confprb{\confpr^{*}}

\def\pn{\pens^{*}}
\def\penInt{\mathfrak{D}}
\def\penH{\mathbb{H}}
\def\pmu{\mathfrak{u}}
\def\Closs{\cc{R}}

\def\dimp{p}
\def\riskb{\riskt_{\rdb}}
\def\dimpp{\dimp+1}
\def\BB{I\!\!B}
\def\vA{\mathtt{v}}

\def\deficiency{\Delta}
\def\spread{\Delta}
\def\dimtotal{\dimp^{*}}

\def\thetav{\bb{\theta}}
\def\thetavs{\thetav^{*}}
\def\thetavc{\thetav'}
\def\thetavd{\thetav^{\circ}}
\def\thetavdc{\thetav^{\sharp}}
\def\dthetavs{\thetav,\thetavs}

\def\thetas{\theta^{*}}
\def\thetac{\theta'}
\def\thetad{\theta^{\circ}}
\def\thetab{\theta^{\dag}}
\def\thetavb{\thetav^{\dag}}

\def\vtheta{\vartheta}
\def\vthetav{\bb{\vtheta}}
\def\prior{\Pi}

\def\Gam{\Xi}
\def\Gam{\mathcal{S}}
\def\RG{R}
\def\Psu{\Upsilon}
\def\Phim{\breve{\Phi}}

\def\Proj{P}

\def\gammavs{\gammav^{*}}
\def\gammavd{\gammav^{\circ}}
\def\etavs{\etav^{*}}
\def\etavd{\etav^{\circ}}
\def\etavc{\etav'}

\def\taus{\tau_{0}}
\def\taup{\lceil \tau \rceil}

\def\sigmas{{\sigma^{*}}}
\def\Sigmas{\Sigma_{0}}

\def\upsilonc{\upsilon'}
\def\upsilond{\upsilon^{\circ}}
\def\upsilonp{{\upsilon}^{*}}
\def\upsilonm{{\upsilon}_{*}}
\def\upsilonvs{\upsilonv^{*}}
\def\upsilons{\upsilon^{*}}
\def\upsilonb{\bar{\upsilon}}
\def\upsilonvd{\upsilonv^{\circ}}

\def\ups{\bb{\upsilon}}
\def\upss{\ups_{0}}
\def\upsc{\ups^{\prime}}
\def\upsd{\ups^{\circ}}
\def\upsdc{\ups^{\sharp}}
\def\upsdu{\ups^{\flat}}

\def\Ups{\varUpsilon}
\def\Upsd{\Ups^{\circ}}
\def\Upss{\Ups_{\circ}}
\def\UpsP{\Ups^{c}}

\def\Thetas{\Theta_{0}}
\def\ThetasGP{\Theta_{0,\GP}}
\def\varthetav{\bb{\vartheta}}

\def\glink{g}

\def\fvs{\fv}
\def\fs{f}
\def\fb{\fv^{\dag}}

\def\uc{\uv'}
\def\ud{\uv^{\circ}}
\def\uvs{\uv^{*}}
\def\us{u^{*}}
\def\vs{v^{*}}

\def\reps{\epsilon}
\def\eps{\epsilon}

\def\repsc{\reps_{0}}
\def\repsb{\reps^{*}}
\def\repsg{g}

\def\lu{\delta}
\def\lub{\bar{\lu}}

\def\Uu{U}
\def\UU{\cc{Y}}
\def\UUM{\cc{M}}
\def\UP{\cc{U}}
\def\up{\mathfrak{u}}

\def\VP{V}
\def\VPc{\VP_{0}}
\def\VPV{\cc{U}}
\def\VPVc{\cc{\VPV}_{0}}
\def\VPGP{\VP_{\GP}}
\def\VPSP{\VP_{\SP}}

\def\VV{H}
\def\GV{\cc{G}}
\def\GVS{S}

\def\VVb{\VV^{*}}
\def\VVc{\VV_{0}}
\def\vv{\bb{h}}
\def\vva{g}
\def\vp{\mathbf{v}}
\def\vpc{\vp_{0}}
\def\VVca{\VV}
\def\Vtt{H}

\def\DG{D}

\def\Vn{V_{0}}
\def\vn{v_{0}}

\def\norm{\mathfrak{c}}
\def\normc{\delta}
\def\norma{c}

\def\egridd{\cc{E}_{\delta}}
\def\penb{\varkappa}

\def\dotzeta{\dot{\zeta}}
\def\mes{\pi}
\def\mesl{\Lambda}
\def\cprr{F}

\def\lambdam{\gm_{1}}
\def\lambdaB{{\lambda}^{*}}
\def\lambdac{{\lambda'}}

\def\cla{{b}}
\def\fis{\mathfrak{a}}
\def\fiss{\fis_{1}}

\def\Vd{{V}}
\def\vd{\bar{v}}

\def\klim{k^{\circ}}
\def\midm{\mid \!}

\def\Ldrift{M}
\def\ldrift{m}
\def\mY{b}
\def\Lvar{D}
\def\lvar{\sigma}

\def\Mubcu{\Upsilon}
\def\Dthetav{\bb{u}}

\def\B{\cc{B}}
\def\BD{\B^{\circ}}
\def\BU{B}
\def\BI{\B^{*}}

\def\mub{\mu^{*}}
\def\mubc{\mu}
\def\mubcb{\mubc^{*}}
\def\Mubc{\mathbb{M}}
\def\Mubcb{\mathrm{M}}

\def\zzc{\zz_{c}}
\def\ww{w}
\def\wwc{\ww_{c}}

\def\norms{\circ} 
\def\rs{\rr_{\norms}}
\def\yys{\yy_{\norms}}
\def\xxs{\xx_{\norms}}
\def\zzs{\zz_{\norms}}
\def\uu{\mathtt{u}}
\def\uus{\uu_{\norms}}
\def\mus{\mu_{\norms}}
\def\gms{\gm_{\norms}}
\def\wws{\ww_{\circ}}

\def\srho{s}
\def\mrho{\varrho}

\def\Lmgf{\mathfrak{M}}
\def\Lmgfb{\Lmgf^{*}}

\def\lmgf{\mathfrak{m}}
\def\lmgfb{\lmgf^{*}}

\def\Expzeta{\mathfrak{N}}
\def\expzeta{\mathfrak{s}}

\def\rr{\mathtt{r}}
\def\rrb{\rr^{*}}
\def\rru{\rr_{\circ}}
\def\rrc{\rr'}
\def\rs{r_{*}}

\def\zz{\mathfrak{z}}
\def\zzb{\tilde{\zz}}
\def\tt{\mathfrak{t}}
\def\zb{z_{\rd}}
\def\zzg{\zz_{1}}
\def\zzQ{\zz_{0}}
\def\zzq{\zz}

\def\Cr{\mathfrak{c}}
\def\Crp{\mathfrak{C}}
\def\Crl{\mathfrak{r}}
\def\Crlp{\mathfrak{R}}
\def\Crlq{\cc{T}}
\def\Crlmu{\cc{M}}

\def\zetah{\zeta_{h}}
\def\GG{G}
\def\HH{H}
\def\pG{p}
\def\pH{q}
\def\hh{H^{*}}

\def\mubch{\mubc_{1}}
\def\rhoh{\rho_{1}}
\def\CoFuh{\CoFu_{1}}
\def\dimh{p_{1}}
\def\VPh{\VP_{1}}
\def\VPt{\VP_{0}}

\def\LLh{L_{1}}
\def\pnndh{\pnnd_{1}}

\def\LCS{C}
\def\Ac{A_{0}}
\def\Ab{A_{\rd}}
\def\DPrb{\DPr_{\rdb}}
\def\DPrm{\DPr_{\rdm}}
\def\Cb{\cc{C}_{\rdb}}
\def\Ub{\cc{U}_{\rdb}}
\def\zetavrb{\zetavr_{\rd}}
\def\xivrb{\breve{\xiv}_{\rd}}
\def\VPrb{\breve{\VP}_{\rdb}}
\def\Larb{\breve{\La}_{\rdb}}
\def\Larm{\breve{\La}_{\rdm}}

\def\deltav{\bb{\delta}}

\def\score{\nabla}
\def\scorer{\breve{\nabla}}

\def\LCS{C}
\def\Ac{A_{0}}
\def\Bc{B_{0}}
\def\AF{A}
\def\Ab{A_{\rdb}}
\def\Am{A_{\rdm}}
\def\DPrc{\DPr_{0}}
\def\DPrb{\DPr_{\rdb}}
\def\DPrm{\DPr_{\rdm}}
\def\Cb{\cc{C}_{\rdb}}
\def\Cm{\cc{C}_{\rdm}}
\def\Ub{\cc{U}_{\rdb}}
\def\deltav{\bb{\delta}}
\def\nuv{\bb{\nu}}
\def\xivrb{\breve{\xiv}_{\rd}}
\def\VPrb{\breve{\VP}_{\rdb}}
\def\Larb{\breve{\La}_{\rdb}}
\def\Lar{\breve{\La}}
\def\Larm{\breve{\La}_{\rdm}}
\def\VH{Q}
\def\VHc{\VH_{0}}
\def\zetavrm{\zetavr_{\rdm}}
\def\N{\mathbb{N}}

\def\Span{\operatorname{span}}
\def\Exc{{\square}}
\def\UUs{U_{\circ}}
\def\errbm{\errb^{*}}
\def\corrDF{\nu}
\def\BBr{\breve{\BB}}
\def\taua{\tau}
\def\AssId{\mathcal{I}}
\def\assId{\iota}
\def\AFD{\cc{A}}

\def\BanX{\cc{X}}
\def\basX{\ev}
\def\apprX{\alpha}
\def\fvs{\fv^{*}}
\def\lkh{\ell}
\def\Bc{B_{0}}
\def\dimn{\dimp_{\nsize}}
\def\betan{\beta_{\nsize}}


\def\xivGP{\xiv_{\GP}}
\def\dimA{\mathtt{p}}
\def\dimAGP{\dimA}
\def\dime{\dimA_{e}}
\def\dimG{\dimA_{\GP}}
\def\dimS{\dimA_{s}}
\def\nubm{\nu_{\rd}}
\def\uub{u_{\rd}}
\def\uubGP{u_{\rd,\GP}}

\def\priorden{\pi}
\def\xivGP{\xiv_{\GP}}
\def\dimAGP{\dimA}
\def\nubm{\nu_{\rd}}
\def\uub{u_{\rd}}
\def\uubGP{u_{\rd,\GP}}

\def\CR{\mathcal{C}}
\def\CRb{\CR_{\rdb}}
\def\vthetavb{\bar{\vthetav}}
\def\Covpost{\mathfrak{S}}

\def\Db{\DP_{+}}
\def\Dm{\DP_{-}}
\def\uvb{\uv_{+}}
\def\uvm{\uv_{-}}
\def\uud{\omega}
\def\taub{\delta}
\def\Lip{L}
\def\Xb{X_{+}}
\def\Xm{X_{-}}
\def\deltam{\delta_{-}}
\def\betauv{\delta}
\def\betab{\betauv_{1}}
\def\betaf{\betauv_{2}}
\def\upsv{\bb{\varkappa}}
\def\upsvb{\bar{\upsv}}
\def\rhob{\varrho}
\def\alpb{\alp_{1}}
\def\betap{\betauv_{3}}
\def\Ec{\E^{\circ}}
\def\ff{f}
\def\fpos{g}
\def\fneg{h}
\def\alpb{\alp_{+}}
\def\alpm{\alp_{-}}

\def\kappak{\kappa}
\def\kappas{\kappak^{*}}
\def\Kappak{\cc{K}}
\def\DPk{\DP_{\kappak}}
\def\VPk{\VP_{\kappak}}

\def\ts{s}
\def\tsv{\bb{\ts}}
\def\mm{\kappa}
\def\mmc{\mm'}
\def\mmd{\mm^{\circ}}
\def\mmo{\mm^{*}}
\def\mmmmo{\mm,\mmo}
\def\mmt{\tilde{\mm}}
\def\mma{\hat{\mm}}
\def\pp{z}

\def\LLL{L_{1}}
\def\LLr{L_{0}}
\def\muL{\mu_{1}}
\def\mur{\mu_{0}}

\def\LmgfL{\Lmgf_{1}}
\def\Lmgfr{\Lmgf_{0}}
\def\Lmgfm{\Lmgf_{1}}

\def\Kappa{\cc{K}}
\def\CoFu{\cc{C}}
\def\CoFuc{\CoFu_{0}}
\def\CoFub{\CoFu^{*}}
\def\CoFuL{\CoFu_{1}}
\def\CoFur{\CoFu_{0}}
\def\CAL{\CA_{1}}
\def\CAr{\CA_{0}}
\def\CAzz{\cc{A}}

\def\pnnL{\pnn_{1}}
\def\pnnr{\pnn_{0}}
\def\ttd{\delta}
\def\alphaL{\alpha_{1}}
\def\alphar{\alpha_{0}}
\def\alpharL{\alpha}
\def\rat{\mathfrak{t}}
\def\mquad{\nquad}
\def\zzL{\zz_{1}}
\def\zzr{\zz_{0}}

\def\mmset{\mathcal{I}}
\def\xex{u}
\def\dcm{q}
\def\dc{g}
\def\dcL{\dc_{1}}
\def\dcr{\dc_{0}}
\def\kk{k}

\def\cpen{\tau}

\def\dens{f}
\def\jj{j}
\def\JJ{\cc{J}}
\def\Zphi{Z}
\def\Zphiv{\bb{\Zphi}}

\def\nuu{\mathfrak{u}}
\def\nud{\mathfrak{u}_{0}}
\def\nun{c_{\nuu}}
\def\rhork{\kullb}
\def\GH{\mbox{GH}}
\def\HYP{\mbox{HYP}}
\def\NIG{\mbox{NIG}}
\def\IR{{\rm I\!R}}
\def\taggr{b}
\def\penm{\boldsymbol{m}}
\def\Crlp{\cc{R}}

\def\Mh{M}
\def\Mht{\Mh^{c}}

\def\Mhh{\Mh^{-}}
\def\Mhc{G}
\def\Lh{L_{1}}
\def\Uh{\cc{U}}
\def\wloc{w}
\def\Bias{B}
\def\bias{b}
\def\ExpzetaU{\Expzeta_{1}}
\def\vpci{\vp_{i,0}}
\def\IFci{\IF_{i,0}}

\def\erqb{\Circle_{\rdb}}
\def\erqm{\Circle_{\rdm}}
\def\errqm{\errm^{*}}
\def\errqb{\errb^{*}}
\def\Nsize{N}
\def\VVD{\VV_{1}}
\def\AA{A}
\def\Wloc{W}

\def\tups{\pen_{0}}
\def\rupd{\rr_{\circ}}
\def\VVb{\VVc}
\def\BP{B}
\def\bp{b}

\def\gps{s}
\def\GK{\cc{G}}

\def\zzGP{\zz_{\GP}}

\def\entrlq{\entrl_{1}}
\def\entrlg{\entrl_{2}}
\def\kb{k^{*}}

\def\rderr{\chi}
\def\Excgr{\diamondsuit}
\def\Excgrb{\diamondsuit^{*}}
\def\Thetat{\bar{\Theta}}
\def\biasGP{\bb{\bias}_{\GP}}
\def\QL{W}
\def\QLG{\mathcal{W}}
\def\BPGP{\QLG_{\GP}}
\def\BBGP{\BB_{\GP}}

\def\xxn{\xx_{\nsize}}
\def\fisGP{\mathtt{w}_{\GP}}
\def\risktGP{\riskt_{\GP}}

\def\dimq{q}
\def\nul{\mathrm{o}}
\def\Thetan{\Theta_{\nul}}
\def\thetavn{\thetav_{\nul}}
\def\thetavsn{\thetavs_{\nul}}
\def\tilden#1{\tilde{#1}_{\nul}}
\def\tildeGP#1{\tilde{#1}_{\GP}}
\def\xivn{\xiv_{\nul}}
\def\xivrGP{\xivr_{\GP}}
\def\DPcc{\DP_{\nul}}
\def\DPnGP{\DP_{1,\GP}}
\def\DPnGPr{\breve{\DP}_{1,\GP}}
\def\nablan{\nabla_{\nul}}
\def\scoren{\score_{\nul}}
\def\AnGP{A_{\nul,\GP}}

\def\testst{T}
\def\TGP{\testst_{\GP}}

\def\VPD{\VP_{2}}

\def\entrlB{\entrl_{1}}
\def\SB{W}
\def\dimq{q}
\def\QQ{\mathbb{H}}
\def\QQg{\QQ_{2}}
\def\QQq{\QQ_{1}}
\def\FF{F}
\def\LaGP{\La_{\GP}}

\def\zzQ{\zz_{\QQ}}
\def\zzAA{\zz_{\FF}}
\def\zzAAA{\zz_{\FF,\SB}}
\def\cdimc{\cdima}
\def\fisGP{\fis_{\GP}}
\def\rdomegab{\rdomega^{*}}

\def\lambdaGP{\lambda_{\GP}}
\def\wGP{\mathrm{w}_{\GP}}

\def\uudm{\mathtt{w}}
\def\lambdaB{\lambda_{\BB}}

\def\lambdav{\bb{\lambda}}
\def\etavd{\etav_{\circ}}
\def\thetavb{\breve{\thetav}}
\def\vthetavd{\Ec \vthetav}
\def\Covd{S_{\circ}}
\def\Covpostd{\Covpost_{\circ}}
\def\IS{\mathcal{I}}
\def\etas{\eta^{*}}
\def\Po{\operatorname{Po}}
\def\IF{\Bbb{F}}
\def\etavb{\bar{\etav}}
\def\etavd{\etav^{\circ}}
\def\Pc{\P^{\circ}}
\def\xxn{\xx_{\nsize}}
\def\CRd{\CR^{\circ}}

\def\CONST{\mathtt{C} \hspace{0.1em}}

\def\dimB{\mathtt{p}_{\BB}}

\def\nub{\nu}
\def\VPD{\VP_{2}}
\def\SB{W}
\def\dimq{q}
\def\dimqb{\dimq^{*}}
\def\QQ{\mathbb{H}}
\def\QQg{\QQ_{2}}
\def\QQq{\QQ_{1}}
\def\FF{F}

\def\qq{z}
\def\qqBB{\qq_{\BB}}
\def\qqQ{\qq_{\QQ}}
\def\qqAA{\qq_{\FF}}
\def\qqAAA{\qq_{\FF,\SB}}
\def\rderr{\chi}
\def\Excgr{\diamondsuit}

\def\Ccb{m}
\def\Ccm{m}
\def\BB{B}

\def\vthetavd{\vthetav^{\circ}}
\def\Indru{\Ind_{\rups}}

\def\BBh{U}
\def\betav{\bb{\beta}}
\def\DD{U}
\def\hsp{\tau}
\def\fiD{a}

\def\Prior{\Pi}
\def\prior{\pi}

\def\In{\mathcal{I}}
\def\KK{\cc{K}}
\def\qqu{\qq^{*}}
\def\ws{\omega}

\newcommand{\tobedone}[1]{\par\textbf{\color{red}To be done:} {\color{magenta}#1}}

\renewcommand{\(}{$\,}
\renewcommand{\)}{\,$}

\def\nquad{\hspace{-1cm}}
\def\eqdef{\stackrel{\operatorname{def}}{=}}
\def\tod{\stackrel{d}{\longrightarrow}}
\def\tow{\stackrel{w}{\longrightarrow}}
\def\toP{\stackrel{\P}{\longrightarrow}}

\newcommand{\cc}[1]{\mathscr{#1}}
\newcommand{\bb}[1]{\boldsymbol{#1}}

\renewcommand{\bar}[1]{\overline{#1}}
\renewcommand{\hat}[1]{\widehat{#1}}
\renewcommand{\tilde}[1]{\widetilde{#1}}

\renewcommand{\Gamma}{\varGamma}
\renewcommand{\Pi}{\varPi}
\renewcommand{\Sigma}{\varSigma}
\renewcommand{\Delta}{\varDelta}
\renewcommand{\Lambda}{\varLambda}
\renewcommand{\Psi}{\varPsi}
\renewcommand{\Phi}{\varPhi}
\renewcommand{\Theta}{\varTheta}
\renewcommand{\Omega}{\varOmega}
\renewcommand{\Xi}{\varXi}
\renewcommand{\Upsilon}{\varUpsilon}
\def\nn{\nonumber \\}

\def\suml{\sum\limits}
\def\supl{\sup\limits}
\def\maxl{\max\limits}
\def\infl{\inf\limits}
\def\intl{\int\limits}
\def\liml{\lim\limits}
\def\Cov{\operatorname{Cov}}
\def\Var{\operatorname{Var}}
\def\arginf{\operatornamewithlimits{arginf}}
\def\argsup{\operatornamewithlimits{argsup}}
\def\argmax{\operatornamewithlimits{argmax}}
\def\argmin{\operatornamewithlimits{argmin}}
\def\val{\operatorname{val}}

\def\D{\boldsymbol{D}}
\def\dd{\operatorname{d}}
\def\tr{\operatorname{tr}}
\def\I{I\!\!I}
\def\R{I\!\!R}
\def\E{I\!\!E}
\def\P{I\!\!P}
\def\X{\mathfrak{X}}
\def\kappa{\varkappa}
\def\Const{\mathrm{Const.} \,}
\def\cdt{\boldsymbol{\cdot}}
\def\tm{\!\times\!}
\def\T{\top}
\def\diag{\operatorname{diag}}
\def\diam{\operatorname{diam}}
\def\rank{\operatorname{rank}}
\def\loc{\operatorname{loc}}

\def\av{\bb{a}}
\def\bv{\bb{b}}
\def\cv{\bb{c}}
\def\dv{\bb{d}}
\def\ev{\bb{e}}
\def\fv{\bb{f}}
\def\gv{\bb{g}}
\def\hv{\bb{h}}
\def\iv{\bb{i}}
\def\jv{\bb{j}}
\def\kv{\bb{k}}
\def\lv{\bb{l}}
\def\mv{\bb{m}}
\def\nv{\bb{n}}
\def\ov{\bb{o}}
\def\pv{\bb{p}}
\def\qv{\bb{q}}
\def\rv{\bb{r}}
\def\sv{\bb{s}}
\def\tv{\bb{t}}
\def\uv{\bb{u}}
\def\vv{\bb{v}}
\def\wv{\bb{w}}
\def\xv{\bb{x}}
\def\yv{\bb{y}}
\def\zv{\bb{z}}

\def\Cv{\bb{C}}
\def\Gv{\bb{G}}
\def\Mv{\bb{M}}
\def\Sv{\bb{S}}
\def\Uv{\bb{U}}
\def\Xv{\bb{X}}
\def\Yv{\bb{Y}}
\def\Zv{\bb{Z}}

\def\alphav{\bb{\alpha}}
\def\epsv{\bb{\varepsilon}}
\def\etav{\bb{\eta}}
\def\gammav{\bb{\gamma}}
\def\varepsilonv{\bb{\varepsilon}}
\def\phiv{\bb{\phi}}
\def\psiv{\bb{\psi}}
\def\tauv{\bb{\tau}}
\def\upsilonv{\bb{\upsilon}}
\def\xiv{\bb{\xi}}
\def\zetav{\bb{\zeta}}

\def\Psiv{\bb{\Psi}}
\def\CONST{\mathtt{C}}

\def\itemv{\vfill\item}
\newenvironment{myslide}[1]
    {\begin{frame}\frametitle{#1}\vfill}
    {\vfill\end{frame}}

\def\vsp{\vspace{0.05\textheight} \vfill}
\def\summarysign{\resizebox{0.08\textwidth}{0.08\textheight}{\includegraphics{summary}}\,}
\def\nix{}
\def\wpu{$\bullet$}

\def\btri{\vfill{\( \blacktriangleright \) }}
\def\btrir{\vfill{\( \blacktriangleright \) }}

\newcommand{\mygraphics}[3]{\begin{center}
    \resizebox{#1\textwidth}{#2\textheight}{\includegraphics{#3}}
    \end{center}
}

\newcommand{\mybox}[3]{\begin{center}
    \resizebox{#1\textwidth}{#2\textheight}{#3}
    \end{center}
}

\newenvironment{eqnh}
{
    \setbeamercolor{postit}{fg=black,bg=hellgelb} 
    \begin{beamercolorbox}[center,wd=\textwidth]{postit} 
    \begin{eqnarray*}}
    {\end{eqnarray*}\end{beamercolorbox}
}

\def\gps{s}
\def\GK{\cc{G}}
\def\Excgr{\diamondsuit}

\def\dimh{m}
\def\LCS{C}
\def\Ac{A_{0}}
\def\Bc{B_{0}}
\def\AF{A}
\def\CF{C}
\def\Ab{A_{\rdb}}
\def\Am{A_{\rdm}}
\def\DPrc{\DPr_{0}}
\def\DPrcp{\DPr_{0,\dimh}}
\def\DPrb{\DPr_{\rdb}}
\def\DPrm{\DPr_{\rdm}}
\def\Cb{\cc{C}_{\rdb}}
\def\Cm{\cc{C}_{\rdm}}
\def\Ub{\cc{U}_{\rdb}}
\def\xivrb{\breve{\xiv}_{\rd}}
\def\VPrb{\breve{\VP}_{\rdb}}
\def\Larb{\breve{\La}_{\rdb}}
\def\Lar{\breve{\La}}
\def\Larm{\breve{\La}_{\rdm}}
\def\VH{Q}
\def\VHc{\VH_{0}}
\def\zetavrm{\zetavr_{\rdm}}

\def\fvh{\bb{\dimh}}
\def\N{\mathbb{N}}
\def\Z{\mathbb{Z}}

\def\iic{\IF}
\def\iif{\breve{\iic}}
\def\dpc{\iic_{\thetav\thetav}}
\def\hhc{\iic_{\etav\etav}}
\def\apc{\iic_{\thetav\etav}}
\def\ifc{\breve{\iic}}

\def\Upsthetav#1{\substack{\\[0.1pt] \upsilonv\in\Ups \\[1pt] \Proj \upsilonv = #1}}
\def\Span{\operatorname{span}}
\def\Exc{{\square}}
\def\UUs{U_{\circ}}
\def\errbm{\errb^{*}}
\def\corrDF{\rho}
\def\BBr{\breve{\BB}}
\def\taua{\tau}
\def\AssId{\mathcal{I}}
\def\AFD{\cc{A}}

\def\BanX{\cc{X}}
\def\basX{\ev}
\def\apprX{\alpha}
\def\fvs{\fv^{*}}
\def\lkh{\ell}
\def\Bc{B_{0}}
\def\h{\frac{1}{2}}
\def\basis{\ev}
\def\Proj{\Pi_{0}}
\def\Projh{\Pi_{\dimh}}

\def\Ij{\mathcal{I}}

\def\NU{\mathbb{H}}

\def\Mn{M_{\nsize}}
\def\bA{\breve{A}}
\def\cA{\bA_{\dimh}}

\def\Sdr{\cc{S}}
\def\xxn{\xx_{\nsize}}

\def\VPr{\breve{\VP}}
\def\BBr{\breve{\BB}}

\def\dimh{m}
\def\LCS{C}
\def\Ac{A_{0}}
\def\Bc{B_{0}}
\def\AF{A}
\def\CF{C}
\def\Ab{A_{\rdb}}
\def\Am{A_{\rdm}}
\def\DPrc{\DPr_{0}}
\def\DPrcp{\DPr_{0,\dimh}}
\def\DPrb{\DPr_{\rdb}}
\def\DPrm{\DPr_{\rdm}}
\def\Cb{\cc{C}_{\rdb}}
\def\Cm{\cc{C}_{\rdm}}
\def\Ub{\cc{U}_{\rdb}}
\def\deltav{\bb{\delta}}
\def\xivrb{\breve{\xiv}_{\rd}}
\def\VPrb{\breve{\VP}_{\rdb}}
\def\Larb{\breve{\La}_{\rdb}}
\def\Lar{\breve{\La}}
\def\Larm{\breve{\La}_{\rdm}}
\def\score{\nabla}
\def\scorer{\breve{\score}}
\def\VH{Q}
\def\VHc{\VH_{0}}
\def\zetavrm{\zetavr_{\rdm}}

\def\fvh{\bb{\dimh}}
\def\N{\mathbb{N}}

\def\IF{{\digamma}}
\def\IF{\Bbb{F}}
\def\iic{\IF}
\def\iif{\breve{\iic}}
\def\dpc{\iic_{\thetav\thetav}}
\def\hhc{\iic_{\etav\etav}}
\def\apc{\iic_{\thetav\etav}}
\def\ifc{\breve{\iic}}

\def\Upsthetav#1{\substack{\upsilonv\in\Ups \\ \Proj \upsilonv = #1}}
\def\Span{\operatorname{span}}
\def\Exc{{\square}}
\def\UUs{U_{\circ}}
\def\errbm{\errb^{*}}
\def\BBr{\breve{\BB}}
\def\taua{\tau}
\def\AssId{\mathcal{I}}
\def\AFD{\cc{A}}

\def\BanX{\cc{X}}
\def\basX{\ev}
\def\apprX{\alpha}
\def\fvs{\fv^{*}}
\def\lkh{\ell}
\def\Bc{B_{0}}
\def\lin{\operatorname{lin}}
\def\h{\frac{1}{2}}
\def\Xv{\bb{X}}
\def\spread{\deficiency}
\def\basis{\ev}
\def\Proj{\Pi_{0}}
\def\Projh{\Pi_{\dimh}}

\def\CONST{\mathtt{C}}
\def\Ij{\mathcal{I}}
\def\etas{\eta^{*}}
\def\nuvs{\nuv^{*}}
\def\nuvc{\nuv'}

\def\omegav{\bb{\phi}}
\def\omegavs{\omegav^{*}}
\def\omegavc{\omegav'}

\def\dimn{\dimp_{\nsize}}
\def\betan{\beta_{\nsize}}
\def\NU{\mathbb{H}}

\def\bA{\breve{A}}
\def\cA{\bA_{\dimh}}

\def\corrDF{\rho}
\def\gps{s}
\def\GK{\cc{G}}
\def\Excgr{\diamondsuit}

\def\dimh{m}
\def\LCS{C}
\def\Ac{A}
\def\Bc{E}
\def\AF{A}
\def\CF{C}
\def\Ab{A_{\rdb}}
\def\Am{A_{\rdm}}
\def\DPc{\DP}
\def\VPc{\VP}
\def\HHc{\HH}
\def\DPrc{\DPr}
\def\DPrp{\DPr_{\dimh}}
\def\DPrb{\DPr_{\rdb}}
\def\DPrm{\DPr_{\rdm}}
\def\Cb{\cc{C}_{\rdb}}
\def\Cm{\cc{C}_{\rdm}}
\def\Ub{\cc{U}_{\rdb}}
\def\xivrb{\breve{\xiv}_{\rd}}
\def\VPrb{\breve{\VP}_{\rdb}}
\def\Larb{\breve{\La}_{\rdb}}
\def\Lar{\breve{\La}}
\def\Larm{\breve{\La}_{\rdm}}

\def\DFc{\DF}
\def\VFc{\VF}

\def\VH{Q}
\def\VHc{\VH}
\def\zetavrm{\zetavr_{\rdm}}

\def\fvh{\bb{\dimh}}
\def\N{\mathbb{N}}
\def\Z{\mathbb{Z}}

\def\iic{\IF}
\def\iif{\breve{\iic}}
\def\DP{{D}}
\def\HH{{H}}
\def\A{{A}}
\def\ifc{\breve{\iic}}

\def\deltar{\delta}

\def\Thetathetav#1{\substack{\\[0.1pt] \upsilonv \in \Ups \\[1pt] \Proj \upsilonv = #1}}
\def\Span{\operatorname{span}}
\def\Exc{{\square}}
\def\UUs{U_{\circ}}
\def\errbm{\errb^{*}}
\def\corrDF{\rho}
\def\BBr{\breve{\BB}}
\def\taua{\tau}
\def\AssId{\mathcal{I}}
\def\AFD{\cc{A}}

\def\BanX{\cc{X}}
\def\basX{\ev}
\def\apprX{\alpha}
\def\fvs{\fv^{*}}
\def\lkh{\ell}
\def\h{\frac{1}{2}}
\def\basis{\ev}
\def\Proj{\Pi_{0}}

\def\Ij{\mathcal{I}}

\def\Mn{M_{\nsize}}
\def\bA{\breve{A}}
\def\cA{\bA_{\dimh}}

\def\Sdr{\cc{S}}
\def\xxn{\xx_{\nsize}}

\def\CONST{\mathtt{C}}
\def\Ij{\mathcal{I}}

\def\etas{\eta^{*}}
\def\zetavs{\zetav^{*}}
\def\zetavc{\zetav'}

\def\omegav{\bb{\phi}}
\def\omegavs{\omegav^{*}}
\def\omegavc{\omegav'}

\def\dimn{\dimp_{\nsize}}
\def\betan{\beta_{\nsize}}

\def\bA{\breve{A}}
\def\cA{\bA_{\dimh}}

\def\corrDF{\rho}
\def\rupf{\rr_{1}}

\def\gmone{\gm_{1}}
\def\rhorb{\rhor_{1}}

\def\upsilonv{\boldsymbol{\upsilon}}
\def\upsilonvs{\boldsymbol{\upsilon}^{*}}
\def\upsilonvd{\boldsymbol{\upsilon}^\circ}
\def\upsilonvc{\upsilonv'}

\def\dimB{\mathtt{p}_{\BB}}

\section{Introduction}
In this work we want to elaborate on the critical dimension in semiparametric profile M-estimation as analyzed in \cite{AASP2013}. Consider a contrast functional \(\LL:\Ups\to \R\)
define for \(\upsilonv=(\thetav,\etav)\)
\begin{EQA}[c]
    \Lr(\thetav) 
    \eqdef 
    \max_{\upsilonv \in \Ups: \, \Proj \upsilonv = \thetav} \LL(\upsilonv).
\label{Lrthetav}
\end{EQA}    
\cite{AASP2013} derive their result for any functional \(\LL\) that satisfies their conditions. The functional need not be the loglikelihood of a parametric family. This will be important for the construction of critical examples in this paper. \\
The object of study, the profile M estimator, is defined as
\begin{EQA}
    \tilde{\thetav}
    &=&
    \argmax_{\thetav \in \Theta} \Lr(\thetav)    
    = 
    \argmax_{\thetav \in \Theta} 
        \max_{\upsilonv \in \Ups: \, \Pi_{\thetav} \upsilonv = \thetav} \LL(\upsilonv).
\label{ttLtt2se}
\end{EQA}    
We define the \emph{semiparametric excess}
\begin{EQA}
    \Lr(\tilde{\thetav}) - \Lr(\thetavs)
    &=&
    \max_{\upsilonv \in \Ups} \LL(\upsilonv) 
    - \max_{\upsilonv \in \Ups: \, \Pi_{\thetav} \upsilonv = \thetavs} \LL(\upsilonv).
\label{excessLr}
\end{EQA}   
The ``target'' value \( \upsilonvs =(\thetavs,\etavs)\)
can defined by 
\begin{EQA}[c]
    \upsilonvs 
    = 
    \argmax_{\upsilonv \in \Ups} \E \LL(\upsilonv) .
\end{EQA}
The key result of \cite{AASP2013} claims that the profile estimator
\( \tilde{\thetav}\) estimates well \( \thetavs \) if the spread \(\Excgr>0\) is small and \(\Lr(\tilde{\thetav}) - \Lr(\thetavs)\approx \|\xivr\|^2\) if \(\sqrt\dimp\Excgr>0\) is small, where \(\|\xivr\|^2>0\) is a quadratic form and \(\dimp\in\N\) the target's dimension. The spread \(\Excgr>0\) is a term that depends on the smoothness of the functional \(\LL\) and exponential bounds of its gradient \(\nabla\LL\). The critical size of the full parameter dimension \(\dimtotal\in\N\) then depends on the exact bounds for \(\Excgr\). 
If the functional is "well behaved" one gets \(\Excgr \asymp \dimtotal/\sqrt{\nsize}\). 
In other words, one needs that ``\({\dimtotal}^{2}/n\) is small'' to obtain 
an accurate non asymptotic version of the Wilks phenomenon and the Fisher Theorem. 

This paper addresses two questions:
1. Are the conditions underlying the result of Theorem \ref{theo: main theo finite dim} necessary for the obtained bound?
2. The error bound for the Fisher expansion differs from that of the Wilks result by a factor of \(\sqrt\dimp\in\N\). Does this difference really exist, i.e. are the results different if \(\dimp= c\dimtotal\)?

We present an example that illustrate that the answer to the first question is "partially yes". If everything is left equal but condition \( {(\breve{\cc{L}}_{0})} \) of \citep{AASP2013} is slightly weakened, then one already needs "\({\dimtotal}^3/n\) is small". Further our second example shows that indeed, once the target dimension is proportional to the full dimension the Wilks result \eqref{eq: Wilks in main theo} becomes substantially harder than the Fisher result \eqref{eq: Fisher in main theo}, i.e. one needs "\({\dimtotal}^3/n\) is small" for the former and ``\({\dimtotal}^{2}/n\) is small'' for the later.

\section{The result of \cite{AASP2013}}
\label{sec: result of AASP}
In this section we summarize the results of \cite{AASP2013}.

Define
\begin{EQA}
\DF\eqdef \nabla^2\E\LL(\upsilonvs),& \VF\eqdef \Cov(\nabla \LL(\upsilonvs), &\zetav(\upsilonv)= \nabla \LL(\upsilonv)-\E\LL(\upsilonv)
\end{EQA}

\cite{AASP2013} prove their main Theorem under the following list of assumptions:

\begin{description}
    \item[\(\bb{(\breve{\LL}_{0})} \)]
    For each \( \rr \le \rups \), 
    there is a constant \( \breve\rddelta(\rr) \) such that
    it holds on the set \( \Upss(\rr) \):
    \begin{EQA}
\label{LmgfquadELGP}
    \|\DP^{-1}\DP^{2}(\upsilonv)\DP^{-1}-I_{\dimp}\|&\le& \breve\rddelta(\rr),\\
     \|\DP^{-1}(\A(\upsilonv)-\A)\HH^{-1}\|&\le& \breve\rddelta(\rr),\\
\left\| \DP^{-1}\A\HH^{-1}\left(I_{\dimh}-\HH^{-1}\HH^2(\upsilonv)\HH^{-1}\right)\right\|
    &\le&\breve \rddelta(\rr).
\end{EQA}

\end{description}

\begin{remark}
This condition describes the local smoothness properties of function \( \E \LL(\upsilonv) \).
We will see that it is necessary for the critical dimension ration "\(\dimtotal/\sqrt{n}\) small".
\end{remark}

\begin{description}
  \item[\( \bb{(\CS \DF)} \)]
    There exist constants \( \nunu>0 \) and \( \gm > 0 \) such that for all 
    \( |\mubc| \le \gm \)
\begin{EQA}[c]
    \sup_{\gammav \in \R^{\dimtotal}} \log\E \exp\left\{ 
        \mubc \frac{\langle \nabla \zeta(\upsilonvd),\gammav \rangle}
                   {\| \VF \gammav \|}
    \right\}
    \le 
    \frac{\nu_{0} ^{2} \mubc^{2}}{2}.
\end{EQA}
\end{description}

\begin{description}
  \item[\( \bb{(\CS \DF_{1})} \)]
    There exists a constant \( \rhor \le 1/2 \), such that for all \( |\mubc| \le \gm \) 
    and all \( 0 < \rr < \rups \)
\begin{EQA}[c]
    \sup_{\upsilonv,\upsilonvc\in\Upss(\rr)}
    \sup_{\|\gammav\|=1} 
    \log \E \exp\left\{ 
         \frac{\mubc \, \gammav^{\T} \DF^{-1} 
         		\bigl\{ \nabla\zetav(\upsilonv)-\nabla\zetav(\upsilonvc) \bigr\}}
         	  {\rhor \, \|\DF (\upsilonv-\upsilonvc)\|}\right\}
    \le 
    \frac{\nu_1^{2} \mubc^{2}}{2}.
\end{EQA}

\end{description}

The global conditions are:

\begin{description}
  \item[\( \bb{(\cc{L}{\rr})} \)] 
     For any \( \rr > \rups\) there exists a value \( \gmi(\rr) > 0 \), 
     such that
\begin{EQA}[c]
    \frac{-\E \LL(\upsilonv,\upsilonvd)}{\|\DF(\upsilonv-\upsilonvd)\|^{2}}
    \ge 
    \gmi(\rr),
    \qquad
    \upsilonv \in \Upss(\rr).
\end{EQA}

  \item[\( \bb{(\CS\rr)} \)] 
    For any \( \rr \ge \rups \) there exists a constant \( \gm(\rr) > 0 \) such that 
\begin{EQA}[c]
    \sup_{\upsilonv \in \Upss(\rr)} \, 
    \sup_{\mubc \le \gm(\rr)} \, 
    \sup_{\gammav \in \R^{\dimtotal}}
    \log\E \exp\left\{ 
        \mubc \frac{\langle \nabla \zeta(\upsilonv),\gammav \rangle}
        {\|\DF\gammav\|}
    \right\}
    \le \frac{\nu_\rr^{2} \mubc^{2}}{2}.
\end{EQA}
\end{description}

Lemma 2.1 of \cite{AASP2013} shows, that these conditions imply the weaker ones  \({(\breve\LL_{0})} \), \({(\breve\CS \DF_{0})} \)  and \({(\breve\CS \DF_{1})} \) that appear bellow:

\begin{lemma}
\label{lem: strong cond imply breve cond} 
Assume \( ({\AssId}) \). Then \( {(\CS \DF_{1})} \) implies \( {(\breve\CS \DF_{1})} \), \( {(\LL_{0})} \) implies \( {(\breve\LL_{0})} \), and \( {(\CS \DF_{0})} \) implies \( {(\breve\CS \DF_{0})} \) with
\begin{EQA}
\breve \gm=\frac{\sqrt{1-\corrDF^2}}{1+\corrDF\sqrt{1+\corrDF^2}}\gm, & \breve \nu=\frac{1+\corrDF\sqrt{1+\corrDF^2}}{\sqrt{1-\corrDF^2}}\nu, & \breve \delta(\rr)=\delta(\rr),\,\text{ and }\breve \omega=\omega.
\end{EQA}
\end{lemma}


\begin{remark}
In this work we concentrate on the conditions \({(\LL_{0})} \), \({(\CS \DF_{0})} \) and \({(\CS \DF_{1})} \) and do not use the refined weaker versions \({(\breve\LL_{0})} \), \({(\breve\CS \DF_{0})} \)  and \({(\breve\CS \DF_{1})} \) as they do not make any difference for the presented examples. 
\end{remark}

We want to cite Theorem 2.1 of \cite{AASP2013}. For this purpose define the \( \dimp \)-vectors \( \scorer_{\thetav} \) and \( \xivr \in \R^{\dimp} \)
\begin{EQA}
	\scorer_{\thetav}
    &=&
    \score_{\thetav} - \Ac \HHc^{-2} \score_{\etav},
    \quad
    \xivr
    \eqdef
    \DPr^{-1} \scorer_{\thetav} ,
\label{scorerdef}
\end{EQA}
and \( \dimp \times \dimp \) matrix \( \DPr^{2}\) as
\begin{EQA}
    \DPr^{2}
    &=&
    \DP^{2} - \Ac \HHc^{-2} \Ac^{\T}. \label{DFse02}
\end{EQA}
The random variable \(\scorer_{\thetav}\in\R^{\dimp} \) is related to the efficient influence function in semiparametric estimation and the matrix \(\DPr^{2}\in \R^{ \dimp \times \dimp}\) equals its covariance in the case of correct specification. 
Define the \emph{semiparametric spread} \(\breve \Excgr(\rr,\xx)>0\) as
\begin{EQA}[c]
\label{eq: def of breve diamond rr}
    \breve\Excgr(\rr,\xx)
    \eqdef
    \left(  \frac{8}{(1-\corrDF^2)^2}\breve\rddelta(\rr) + \, 6 \nu_{1}  \breve\omega \zzq(\xx,2\dimtotal+2\dimp) \right)\rr,
\label{Exceqrrrhdef}
\end{EQA}
where \( \breve\rddelta(\rr) \) is shown in the condition \({ (\breve\LL_{0})} \), 
the constants \(  \breve\rhor \), \(\nu_1 \) are from condition \({ (\breve\CS \DF_{1})} \). 
The value \( \zzq(\xx,2\dimtotal+2\dimp) \) is related to the entropy of the unit ball
in a \( \R^{\dimtotal+\dimp} \)-dimensional Euclidean space
\begin{EQA}[c]
    \zzq(\xx,2\dimtotal+2\dimp)
    \eqdef
    \begin{cases}
    \sqrt{2 (\xx +2\dimtotal+2\dimp)} 
        & \text{if } \sqrt{2 (\xx + 2\dimtotal+2\dimp)} \le \gm, \\
    \gm^{-1} (\xx + 2\dimtotal+2\dimp) + \gm/2 
        & \text{otherwise} ,
  \end{cases}
\label{eq: def of entropy term in main chapter}
\end{EQA}  
and one can apply \(\zzq(\xx,\dimtotal)\cong \sqrt{\xx+\dimtotal}\)
 for moderate choice of \(\xx>0\); see Appendix C of \cite{AASP2013}. 
The value \(\breve \Excgr(\rr,\xx) \) measures the quality of a linear approximation to \(\scorer\LL(\upsilonv)-\scorer\LL(\upsilonvs)\) in the local vicinity the local vicinity 
\( \Upss(\rr) = \bigl\{ \upsilonv\in\Ups \colon \|\DFc(\upsilonv-\upsilonvd)\|\le \rr \bigr\} \). 
Finally we introduce \( \tilde{\upsilonv}_{\thetavs}\in\Ups \), which maximizes \( \LL(\upsilonv,\upsilonvs) \) subject to
\(\Proj \upsilonv = \thetavs \):
\begin{EQA}[c]
\label{tuthLLuus}
    \tilde{\upsilonv}_{\thetavs} \eqdef (\thetavs,\tilde \etav_{\thetavs})
    \eqdef 
    \argmax_{\Thetathetav{\thetavs}} 
    \LL(\upsilonv,\upsilonvs),
\end{EQA}
and define the radius \(\rups>0\)
\begin{EQA}[c]\label{eq: def of rups}
\rups(\xx)\eqdef \inf_{\rr>0}\left\{\P(\tilde{\upsilonv},\tilde{\upsilonv}_{\thetavs}\in\Upss(\rr))\ge 1-\ex^{-\xx}\right\},
\end{EQA}
which we set to infinity if \(\tilde{\upsilonv}=\{\,\}\) or \(\tilde{\upsilonv}_{\thetavs}=\{\,\}\).

\begin{theorem}
\label{theo: main theo finite dim}
Assume 
\({(\breve\CS \DF_{1})} \), \({(\breve\LL_{0})}\), and \({(\AssId)}\) with 
a central point \(\upsilonvd =\upsilonvs\) 
and some matrix \( \DFc^{2} \). Further assume that the sets of maximizers \(\tilde\ups,\,\tilde{\upsilonv}_{\thetavs}\) are not empty.
Then it holds on a set \(\Omega(\xx)\subseteq\Omega \) of probability at least
\( 1-2\ex^{-\xx} \) for the profile MLE \( \tilde{\thetav} \)
\begin{EQA}
	\bigl\| 
        \DPr \bigl( \tilde{\thetav} - \thetavs \bigr) 
        - \xivr 
    \bigr\|
    &\le& 
    \breve\Excgr(\rups,\xx) ,
\label{eq: Fisher in main theo}
	\\
    \bigl| 2 \Lr(\tilde{\thetav},\thetavs) - \| \xivr \|^{2} \bigr|
    &\le&
     8\left(\|\xivr\|+\breve\Excgr(\rups,\xx)\right)\breve\Excgr(2(1+\corrDF)\rups,\xx)+ \breve\Excgr(\rups,\xx)^2,
\label{eq: Wilks in main theo}
\end{EQA}
where the spread \(\breve \Excgr(\rups,\xx) \) is defined in \eqref{eq: def of breve diamond rr}.
\end{theorem}

The critical size of \(\dimtotal\) depends on the exact bounds on \(\delta(\cdot),\rhor\). 
If \(\breve\delta(\rr)/\rr \asymp \breve\rhor \asymp 1/\sqrt{\nsize}\) one gets \(\Excgr(\rr,\xx) \asymp \dimtotal/\sqrt{\nsize}\). \cite{AASP2013} discuss the critical dimension in this setting and in our example we will have \(\breve\delta(\rr)/\rr \asymp \breve\rhor \asymp 1/\sqrt{\nsize}\) as well.

\subsection{Critical smoothness}
This section addresses the necessary smoothness to ensure the bound \({\dimtotal}^2/n\) for Theorem \ref{theo: main theo finite dim}. We show that the following slightly weaker version of \( \bb{(\LL_{0})} \) already allows to find examples that satisfy all conditions of Section  2.1 of \citep{AASP2013} but for which the critical ratio is \({\dimtotal}^3/n\to 0\).

Consider:
\begin{description}
  \item[\( \bb{(\breve{\cc{L}}_{0})'} \)]
    There exists a symmetric \( \dimtotal \times \dimtotal \)-matrix \( \DF^{2} \) such 
    that 
    such that it holds on 
    the set \( \Upss(\rups) \) for all \( \rr \le \rups \)
\begin{EQA}[c]
    \left\lvert 
        \frac{ \E \LL(\ups,\upss) - \|\DF(\upsilonv - \upsilonvs)\|^2}
             {\|\DF(\upsilonv - \upsilonvs)\|^2}
    \right\rvert
    \le 
    \rddelta(\rr).
\end{EQA}
\end{description}
and note that \( \bb{(\LL_{0})} \) implies \(\bb{(\cc{L}_{0})'}\) but not the other way round.

\begin{remark}
This condition appears in \cite{SP2011} and allows to prove Theorem \ref{theo: main theo finite dim} from above with an error bound \(\triangle(\rr,\xx)=O({\dimtotal}^{3/2}/\sqrt{n})\) instead of \(\Excgr(\rr,\xx)\).
\end{remark}

We write \( \dimtotal = \dimn \).
We show that the condition \({\dimtotal}^3/n\to 0\) is critical for the class of models
satisfying the conditions of Section~\ref{sec: result of AASP} with \( \bb{(\LL_{0})} \) replaced by \( \bb{(\LL_{0})'} \). 
%
Namely, we present an example 
in which the behavior of the profile MLE \( \tilde{\thetav} \) heavily depends on the 
value \( \betan = \sqrt{\dimn^{3} / \nsize} \ge \beta > 0 \). 
If \( \betan \to 0 \), then the conditions of Section~\ref{sec: result of AASP} 
are satisfied yielding asymptotic efficiency of \( \tilde{\thetav} \). 
At the same time, if \( \betan \ge \beta > 0 \), then
the MLE \( \tilde{\thetav} \) is not anymore root-n consistent.

Assume that \(\dimn/\sqrt n\to 0\).
Let a random vector \( \Xv \in \R^{\dimn} \) follow 
\( \Xv \sim \ND(\upsilonvs,\nsize^{-1} \Id_{\dimn}) \).
Take for simplicity \( \upsilonvs = 0 \) and let \( \P = \P_{0} \) denote the distribution 
of \( \Xv \).
Introduce a special set \( \Sdr \subset \R^{\dimn} \) with
\begin{EQA}
    \Sdr
    & \eqdef &
    \left\{ \upsilonv=(\upsilon_1,\ldots,\upsilon_{\dimn}): \, \upsilon_{1} = \frac{z}{2} \sqrt{\betan/\nsize} , \, z\in\Z \right\}\\
    &&\cap \,\Upss\left(\sqrt{2 \dimn/\nsize}+\frac{1}{2} \sqrt{\betan/\nsize}\right).
\label{Sdrsetse}
\end{EQA}
We denote by \(\Sdr_{\delta} \) its \( \delta \)-vicinity:
\begin{EQA}
    \Sdr_{\delta}
    & \eqdef &
    \{\upsilonv:\,\dist(\upsilonv,\Sdr) < \delta \} ,
\label{Sdrdelta}
\end{EQA}    
where \( \dist(\upsilonv,\Sdr) \) is the Euclidean distance from the point 
\( \upsilonv \) to the set \( \Sdr \).
Also \( \Sdr_{\delta}^{c} \) stands for the complement of \( \Sdr_{\delta} \).
Below we fix \( \delta = 1/\nsize \).
Consider a special parametric quasi log-likelihood ratio 
\( \LL(\upsilonv,0) \) defined as
\begin{EQA}
    \LL(\upsilonv,0)
    &=& 
    \nsize \Xv^{\T} \upsilonv - \nsize \|\upsilonv\|^{2}/2 
        + \nsize f(\upsilonv) \|\upsilonv\|^{3} .
\end{EQA}
Here \( f: \R \mapsto \R \) is a smooth function with
\begin{EQA}[c]
    f(\upsilonv)
    =
    \begin{cases}
      1 & \upsilonv \in \Sdr, \\
      0 & \upsilonv \in \Sdr_{\delta}^{c}.
    \end{cases}
\end{EQA}
Below we consider the problem of estimating the first component 
\( \theta \eqdef \upsilon_{1} \in \R \).
Since by assumption \(\dimn/\sqrt{\nsize}\to 0\) it holds for \( \nsize \) large enough and for any \( \upsilonv \) with 
\( \| \upsilonv \|^{2} \le 4\dimn/\nsize+\betan/\nsize \) that
\( \nsize \|\upsilonv\|^{2}/2 \ge \nsize f(\upsilonv) \|\upsilonv\|^{3} \) and thus
\begin{EQA}[c]
    \argmax_{\upsilonv} \E \LL(\upsilonv)
    =
    \argmin_{\upsilonv} 
    \bigl\{ \nsize \|\upsilonv\|^{2}/2 - \nsize f(\upsilonv) \|\upsilonv\|^{3} \bigr\}
    =
    0.
\end{EQA}
It is easy to see that all conditions from Section~\ref{sec: result of AASP} except \( \bb{(\cc{L}_{0})} \) are satisfied 
with \( \rhor \cong 1/\sqrt{n} \) and
\begin{EQA}
   \DF^{2} 
    = 
    \VF^{2}
    = 
    \nsize \Id_{\dimn}, & & \Upss(\rr)=\{\|\ups\|\le \rr/\sqrt{n}\} .
\end{EQA}
But clearly \( \bb{(\cc{L}_{0})'} \) is met with \(\delta(\rr)=\rr/\sqrt{n}\).
It is straightforward to see that 
\begin{EQA}[c]
    \DPrc
    =
    \sqrt{\nsize},
    \qquad \scorer(\LL-\E\LL)
    = 
    \score_{\thetav} (\LL-\E\LL)
    =
    \nsize X_{1}, 
    \,\text{ and } 
    \xivr = \sqrt{\nsize} X_{1}.
\end{EQA}
The next result shows that in this example the critical ratio reads \( \betan = \sqrt{\dimn^{3}/\nsize} \), i.e. iff it is not small,
the profile MLE \( \tilde{\theta} \) is not root-\( \nsize \) consistent.

\begin{theorem}
\label{Tcritexam}
If \( \betan^{2} = \dimn^{3} / \nsize \to 0 \) then 
\begin{EQA}[c]
    \| \DPrc(\tilde{\theta} - \thetas) - \xivr \|
    = 
    \sqrt{\nsize} |\tilde \upsilon_{1} - X_{1} |
    \to 0.
\end{EQA}
Suppose that \( \betan \to (6c)^{2} \) for some \( c > 0 \).
Let also \( \nsize \) be large enough to ensure
\begin{EQA}[c]
    \frac{2^{1/3}-1}{2^{1/6}}\sqrt{\dimn/\nsize}
    \ge  
    \frac{1}{2} \left({\dimn}/{\nsize}\right)^{3/4} .
\end{EQA}
There exists a positive \( \alpha > 0 \) such that 
it holds with a probability exceeding \( \alpha \) 
\begin{EQA}[c]
    \| \DPrc(\tilde{\theta} - \thetas) - \xivr \|
    \ge 
    \frac{1}{6} \betan^{1/2} - \frac{1}{\sqrt{\nsize}}
    \ge c - o_{\nsize}(1).
\end{EQA}
If \(\betan \to \infty\), then 
\begin{EQA}[c]
    \| \DPrc(\tilde{\theta} - \thetas) - \xivr \|
    \toP 
    + \infty,
\end{EQA}
where \( \toP \) means convergence in probability.
\end{theorem}

In short: we have shown that - everything else left unchanged - a smoothness conditions of the kind of \( \bb{(\cc{L}_{0})} \), i.e. qualified smoothness of second derivatives, is necessary to ensure that "\(\dimtotal/\sqrt{n}\) is small" suffices to get accurate results in Theorem \ref{theo: main theo finite dim}.

\subsection{Difference between Wilks and Fisher}
\label{sec: example}
This section discusses the issue of \emph{critical dimensions} if the target dimension \( \dimp= c\dimtotal \) for some \(c>0\). We again write \( \dimtotal = \dimn \). 
In this case Theorem~\ref{theo: main theo finite dim} - and assuming that \(\delta(\rr)/\rr\cong \omega\cong 1/\sqrt n\) - requires that \( \dimn = o(\nsize^{1/3}) \) or \( \dimn = o(\nsize^{1/2})\) to obtain non asymptotic versions of the Wilks phenomenon and the Fisher Theorem respectively.
Here we show that this difference actually occurs on the class of models
satisfying the conditions of Section~\ref{sec: result of AASP}.

%
Our example shows critical behavior in the sense that, when \(\dimn^{3}/ n\nrightarrow 0\) we find for each
\(n\in\N\) large enough a set \(\mathcal A\subset \Omega\) of positive probability on which the profile log likelihood ratio does not
converge to a chi-square random variable. In accordance with the results of Theorem \ref{theo: main theo finite dim} the estimator is efficient if
\(\dimn^{2}/ n\rightarrow 0\) and the Wilks phenomenon occurs if \(\dimn^{3}/n\rightarrow 0\).

Assume \(\dimn=2\dimh\) and take as target \(\thetav:=\Pi_{2}\upsilonv\in \R^{\dimh}\), where \(\Pi_{2}: \R^{\dimn}\to \R^{\dimh}\) denotes the orthogonal projection on the first \(\dimh\in\N\) components.
Further assume that \(\dimn^2/n\rightarrow 0\). We use a
miss specified model, i.e. we take standard normal observations on \(\R^{\dimn}\) but assume that the ME is derived from the
correct loglikelihood function altered by an additional term. Consider
\begin{EQA}[c]
\LL(\upsilonv,0)= n\mathbf X^\T \upsilonv-n\|\upsilonv\|^2/2+ f(\upsilonv)n\|\upsilonv\|^{3}/3,
\end{EQA}
where
\begin{EQA}[c]
\Xv \sim N\left(0,\frac{1}{n}I_{\dimn}\right),
\end{EQA}
where \(f:\R^{\dimn}\mapsto \R\) is some smooth function with
\begin{EQA}[c]
f(\upsilonv)=\begin{cases}
      0 & \upsilonv_1=0,\\
      1 & \mathcal S:=\left\{\|\Pi_{2}\upsilonv\| \ge \frac{2}{L}\sqrt{\frac{\dimn}{n}}\right\}\cap B_{2\sqrt{\frac{\dimn}{n}}}(0),
     \end{cases}
\end{EQA}
where \(L>0\). More precisely we set for any \(\upsilonvd\in\R^{\dimn}\)
\begin{EQA}[c]
\label{eq: def of f}
f(\upsilonvd)=\varphi_{\left\{\|\Pi_{2}\upsilonv\| \ge \frac{2}{L}\sqrt{\frac{\dimn}{n}}\right\}}(\upsilonvd)1_{B_{2\sqrt{\frac{\dimn}{n}}}(0)}(\upsilonvd),
\end{EQA}
where 
\begin{EQA}[c]
\varphi_{\left\{\|\Pi_{2}\upsilonv\| \ge \frac{2}{L}\sqrt{\frac{\dimn}{n}}\right\}}(\upsilonvd)=\int_{\R} 1_{\{\|\Pi_{2}\upsilonv\| \le \frac{1}{L}\sqrt{\frac{\dimn}{n}}\}}(\upsilonv_1)K_{\frac{1}{L}\sqrt{\frac{\dimn}{n}}}(\upsilonvd_1-\upsilonv_1)d\upsilonv_1\,
\end{EQA}
where \(K\) is a smooth kernel with support on \([-1,1]\) and 
\begin{EQA}[c]
K_{h}(x):=\frac{1}{h}K\left(\frac{x}{h}\right).
\end{EQA}


\begin{theorem}
\label{the: critical dim wilks}
In the above model the conditions of Section \ref{sec: result of AASP} are satisfied yielding \(\Excgr(\rups,\xx)=o( \dimn/ \sqrt{n}) \). The Fisher theorem holds true if \(\dimn^{2}/ n\rightarrow 0\). Further the Wilks phenomenon occurs iff \(\dimn^{3}/ n\rightarrow 0\).
\end{theorem}

\subsection{Proof of Theorem~\ref{Tcritexam}}
We only sketch the proof of the first claim as it is rather uninteresting. Note that 
\begin{EQA}[c]
\P(n\|\Xv\|^2\ge 4\dimn)\to 0,
\end{EQA}
which implies that \(\tilde\upsilonv\in \Upss(2\sqrt{\dimn})\). On \(\Upss(2\sqrt{\dimn})\)
\begin{EQA}[c]
\label{eq: quad bracketing}
n\upsilonv^\T\Xv-n\|\upsilonv\|^2/2-\sqrt{\betan}\le \LL(\upsilonv)\le n\upsilonv^\T\Xv-n\|\upsilonv\|^2/2+\sqrt{\betan}.
\end{EQA}
Maximizing on the left hand side of \eqref{eq: quad bracketing} and plugging in \(\tilde\upsilonv\) on the right hand side we get
\begin{EQA}[c]
\|\DF(\tilde\upsilonv-X)\|^2/2=n\|Xv\|^2/2-n\tilde\upsilonv^\T\Xv+n\|\tilde\upsilonv\|^2/2\le 2\sqrt{\betan}.
\end{EQA}
This gives the claim:
\begin{EQA}
\| \DPrc(\tilde{\thetav} - \thetavs) - \xivr \| ^2
    &\le& \|\DF(\tilde\upsilonv-X)\|^2\le 2\sqrt{\betan}\to 0.
\end{EQA}

For the other claims we first show that for \( \nsize \) 
large enough, the MLE 
\( \tilde{\upsilonv} \in \R^{\dimn} \) belongs with probability close to one 
to the \( \delta = 1/\nsize \) vicinity \( \Sdr_{\delta} \) of the set \( \Sdr \) 
from \eqref{Sdrsetse}.
The second step is to show that with a probability exceeding a fixed 
constant \( \alpha > 0 \),
the profile MLE \( \tilde{\theta} \) differs significantly from \( X_{1} \) which is the 
profile MLE in the linear Gaussian model. 
The third step focuses on the case \( \betan \to \infty \). 

1. First we show that for \( \nsize \) large enough, 
the MLE \( \tilde{\upsilonv}\in\R^{\dimn} \) lies in \( \Sdr_{\delta} \) 
with probability close to one. 
For this we check that the maximum of \( \LL(\upsilonv) \) on 
\( \Sdr_{\delta}^{c} \) is smaller than a similar maximum on \( \Sdr \)
for ``typical'' values of \( \Xv \) and \( \nsize \) large enough.
Indeed, for any point \( \upsilonv \in \Sdr_{\delta}^{c} \) 
\begin{EQA}
    \LL(\upsilonv,0)
    & \le & 
    \max_{\upsilonv \in \Sdr_{\delta}^{c}} \LL(\upsilonv,0)
    =
    \max_{\upsilonv \in \Sdr_{\delta}^{c}}
        \bigl\{ \nsize \Xv^{\T} \upsilonv - \nsize \|\upsilonv\|^{2}/2 \bigr\}
    \\
    & \le & 
    \max_{\upsilonv \in \R^{\dimn}} 
        \bigl\{ \nsize\Xv^{\T} \upsilonv-\nsize\|\upsilonv\|^{2}/2 \bigr\}
    = 
    \frac{\nsize}{2} \| \Xv \|^{2}.
\end{EQA}  
Further, introduce a random set of ``typical'' values \( \Xv \):
\begin{EQA}[c]
    \LCS_{1} 
    \eqdef 
    \left\{\Xv: \,\, 
        \frac{1}{2} \left(\frac{\dimn}{\nsize}\right)^{3/2} < \| \Xv \|^{3} 
        < \left(\frac{2 \dimn}{\nsize}\right)^{3/2}, 
        \text{ and } |X_{1}| \le 1
    \right\} .
\end{EQA}
It is straightforward to see that 
\( \P\bigl( \Xv \in \LCS_{1} \bigr) \) is exponentially close to one for \( \nsize \)
large. 
Below we assume that \( \Xv \in \LCS_{1} \) and 
study the value \( \LL(\upsilonv,0) \) for \( \upsilonv \in \Sdr \).
Let also \( \nsize \) be large enough to ensure that
\begin{EQA}[c]
\label{eqconditionn}
    \frac{2^{1/3}-1}{2^{1/6}}\left(\frac{\dimn}{\nsize}\right)^{1/2} 
    \ge  
    \frac{1}{2}\left(\frac{\dimn}{\nsize}\right)^{3/4} 
    = 
    \frac{1}{2} \sqrt{\betan/\nsize} .
\end{EQA}
Introduce \( \Xv_{\Sdr} \) as the closest point in 
\( \Sdr \) to \( \Xv \)  with \( |\upsilon_{1}| \ge |X_{1}| \). 
This point always exists by the definition of \( \Sdr \). 
Denote 
\begin{EQA}[c]
    \delta(\Xv) = \| \Xv - \Xv_{\Sdr} \| 
    = 
    |X_{1} - \upsilon_{1}| .
\label{deltaXv1}
\end{EQA}    
By construction of \( \Sdr \), it holds 
\( \delta(\Xv) \le 0.5 \sqrt{\betan/\nsize} \) for \( \Xv \in \LCS_{1} \). 
For \( \nsize \) satisfying \eqref{eqconditionn} this also yields
\( \bigl[ \| \Xv \| - \delta(\Xv) \bigr]^{3} \ge 1/2 \| \Xv \|^{3} \).
Now we have for \( \Xv \in \LCS_{1} \)
\begin{EQA}
    \max_{\upsilonv \in \Sdr} \LL(\upsilonv,0)
    &\ge&  
    \LL(\Xv_{\Sdr},0)
    \\
    &\ge& 
    \nsize\| \Xv \|^{2} - \nsize |X_{1}| \delta(\Xv)
    - \frac{\nsize}{2} \bigl\{ \| \Xv \|^{2} - 2 |X_{1}| \delta(\Xv) + \delta^{2}(\Xv) \bigr\}
    \\
    &&
    \qquad + \, 
    \nsize \bigl\{ \| \Xv \|^{2} - 2 |X_{1}| \delta(\Xv) + \delta^{2}(\Xv) \bigr\}^{3/2}
    \\
    & \ge &
    \frac{\nsize}{2} \| \Xv \|^{2} - \nsize \delta^{2}(\Xv)
    + \nsize \bigl\{ \| \Xv \| - \delta(\Xv) \bigr\}^{3}
    \\
    &>&
    \frac{\nsize}{2} \| \Xv \|^{2} - \frac{\betan}{4} + \frac{\nsize}{2} \| \Xv \|^{3}
    >
    \frac{\nsize}{2} \| \Xv \|^{2}
    =
    \max_{\upsilonv \in \Sdr_{\delta}^{c}} \LL(\upsilonv,0).
\end{EQA}
This implies \( \tilde{\upsilonv} \in \Sdr_{\delta} \).

2. Now we discuss the case when \( \betan^{2} = \dimn^{3} / \nsize \to (6c)^{4} \) for some 
\( c \ge 0 \) and show that 
the profile MLE \( \tilde{\thetav} \) 
deviates significantly from \( X_{1} \)
on a random set of positive probability. 
Define for each \( \nsize \in \N \) 
\begin{EQA}[c]
    \LCS_{\nsize} 
    \eqdef 
    \LCS_{1} \cap \left\{ \| \Xv - \Xv_{\Sdr} \|
    \ge \frac{1}{6} \sqrt{\betan/\nsize} \right\}
    =
    \LCS_{1} \cap \left\{ |X_{1} - X_{\Sdr,1}|
    \ge 
    \frac{1}{6} \sqrt{\betan/\nsize} \right\} .
\end{EQA}
It is easy to see that \( \P(\LCS_{\nsize}) \ge \alpha \)  
for some fixed \( \alpha > 0 \) and all \( \nsize \).
It remains to note that on the set \( \LCS_{\nsize} \) it holds under 
\eqref{eqconditionn} 
\begin{EQA}
    \| \DPrc(\tilde{\thetav} - \thetavs) - \xivr \| 
    &=&  
    \sqrt{\nsize} |\tilde{\upsilon}_{1} - X_{1}| 
    \\
    &\ge& 
    \sqrt{\nsize} |X_{1} - X_{\Sdr,1}| - \sqrt{\nsize} / \nsize
    \\
    & \ge &
    \frac{1}{6} \betan^{1/2} - \frac{1}{\sqrt{\nsize}}
    \to 
	\begin{cases}
	    \infty & {\dimn^{3}}/{\nsize} \to \infty,
        \\
	    c & {\dimn^{3}}/{\nsize} \to (6c)^{4}, 
    \end{cases} 
\end{EQA}
which yields the claim.

3. Finally consider the case when \( \betan \to \infty \).
Fix any sequence \( c_{\nsize} \) such that 
\( c_{\nsize} \to 0 \) and \( c_{\nsize} \betan \to \infty \),
e.g. \( c_{\nsize} = \betan^{-1/2} \).
Consider the random set 
\begin{EQA}[c]
    \LCS_{\nsize} 
    \eqdef 
    \LCS_{1} \cap \left\{ \| \Xv - \Xv_{\Sdr} \|
    \ge \frac{c_{\nsize}}{6} \sqrt{\betan/\nsize} \right\}
    =
    \LCS_{1} \cap \left\{ |X_{1} - X_{\Sdr,1}|
    \ge 
    \frac{c_{\nsize}}{6} \sqrt{\betan/\nsize} \right\}.
\end{EQA}
Then \( \P(\LCS_{\nsize}) \to 1 \) and on \( \LCS_{\nsize} \) 
\begin{EQA}
    \| \DPrc(\tilde{\theta} - \thetas) - \xivr \|
    & \ge & 
    \frac{c_{\nsize}}{6} \betan^{{1}/{2}}-\frac{1}{\sqrt{\nsize}}
    \to \infty,
\end{EQA}
as required.

\subsection{Proof of Theorem~\ref{the: critical dim wilks}}

Since by assumption \(\dimn^2/\nsize\to 0\) and the support of \(f(\upsilonv)\) is contained in \(B_{2\sqrt{\dimn}/\sqrt{n}}(0)\) by the choice of \(\hat K\) it holds for \( \nsize \) large enough and for any \( \upsilonv \) with 
\( \| \upsilonv \|^{2} \le 4\dimn/\nsize \) that
\( \nsize \|\upsilonv\|^{2}/2 \ge \nsize f(\upsilonv) \|\upsilonv\|^{3} \) and thus
\begin{EQA}[c]
    \argmax_{\upsilonv} \E \LL(\upsilonv)
    =
    \argmin_{\upsilonv} 
    \bigl\{ \nsize \|\upsilonv\|^{2}/2 - \nsize f(\upsilonv) \|\upsilonv\|^{3}/3 \bigr\}
    =
    0.
\end{EQA}

Appart from \(\bb{(\LL_{0})}\) it is easy to see that all conditions are satisfied with \(\gmi=1\) and \(\rddelta(\rr)/\rr\cong \omega \cong 1/\sqrt{n}\) if we set
\begin{EQA}[c]
\DF^2=\VF^2=nI_{\dimn}.
\end{EQA}
It is straightforward to see that 
\begin{EQA}[c]
    \DPrc
    =
    \sqrt{\nsize},
    \qquad \scorer(\LL-\E\LL)
    = 
    \score_{\thetav} (\LL-\E\LL)
    =
    \nsize X_{1},
    \,\text{ and } 
    \xivr = \sqrt{\nsize} X_{1}.
\end{EQA}
Consequently Theorem \ref{theo: main theo finite dim} gives efficiency of the profile if \(\frac{\dimn^{2}}{ n}\rightarrow 0\) and the Wilks phenomenon if \(\frac{\dimn^{3}}{ n}\rightarrow 0\). In the following we will first show that condition \(\bb{(\LL_{0})}\) is satisfied, then that the Fisher theorem holds if \(\frac{\dimn^{2}}{ n}\rightarrow 0\) and finally that \(\frac{\dimn^{3}}{ n}\rightarrow 0\) is indeed necessary to obtain the Wilks phenomenon.

\subsubsection{Condition \(\bb{(\LL_{0})}\)}
We will show that \(\nabla^2\E\LL\) is Lipshitz continuous on \(\Upss(\rups)=B_{2\sqrt{\dimn/n}}(0)\) with Lipshitz constant \(n\tilde L>0\) where \(\tilde L\) is independent of \(n,\dimn\). This gives \(\bb{(\LL_{0})}\) with \(\rddelta(\rr)=\tilde L\rr/\sqrt n\). For this purpose it suffices to consider the Lipshitz continuity of \(\nabla^2g(\upsilonv):=\nabla^2( f(\upsilonv)n\|\upsilonv\|^{3})\). We neglect the indicator \(1_{B_{2\sqrt{\dimn/n}}(0)}(\cdot)\) as we only have to consider smoothness on \(\Upss(\rups) \). We have for two points \(\upsilonv,\upsilonvd\in\Upss\)
\begin{EQA}
  \frac{1}{n}\|\nabla^2g(\upsilonv)-\nabla^2g(\upsilonvd)\|&\le&\|\nabla^2 f(\upsilonv)\|\upsilonv\|^{3}-\nabla^2f(\upsilonvd)\|\upsilonvd\|^{3} \|\\
   &&+\|\nabla f(\upsilonv)\|\upsilonv\|\upsilonv^\T- \nabla f(\upsilonvd)\|\upsilonvd\|{\upsilonvd}^\T\|\\
   &&+\|\frac{f(\upsilonv)}{\|\upsilonv\|}\upsilonv\upsilonv^\T-\frac{f(\upsilonvd)}{\|\upsilonvd\|}\upsilonvd{\upsilonvd}^\T \|.
\end{EQA}
Denote by \(L_{\|\cdot\|^3|_{\Upss}}\) the Lipshitz constant of \(\|\cdot\|^3\) restricted to \(\Upss(\rups)\), which is independent of \(n,\dimn\in\N\) because the set \(\Upss(\rups)\subset B_1(0)\) for \(n\in\N\) large enough.
We estimate 
\begin{EQA}
&&\nquad\|\nabla^2 f(\upsilonv)\|\upsilonv\|^{3}-\nabla^2f(\upsilonvd)\|\upsilonvd\|^{3} \|\\
&\le& \|\nabla^2 f(\upsilonv)-\nabla^2f(\upsilonvd)\|\|\upsilonv\|^{3}+\|\nabla^2 f(\upsilonvd)\|\|\|\upsilonv\|^{3}-\|\upsilonvd\|^{3} \|\\
  &\le&8\left(\frac{\dimn}{n}\right)^{3/2}\|\nabla^3f\|_{\infty}\|\upsilonv-\upsilonvd\|+\|\nabla^2 f\|_{\infty}L_{\|\cdot\|^3|_{\Upss}}\|\upsilonv-\upsilonvd\|.
\end{EQA}
By the definition \eqref{eq: def of f} we find that
\begin{EQA}[c]
\|\nabla^3f\|_{\infty}\le L^3\left(\frac{n}{\dimn}\right)^{3/2} \left\|\int_{\R} K^{(3)}(\upsilonv)d\upsilonv \right\|=:\CONST\left(\frac{n}{\dimn}\right)^{3/2},
\end{EQA}
with a constant \(\CONST\in\R\) that does not depend on \(n,\dimn\in\N\).
With the same arguments for the other terms we find
\begin{EQA}[c]
\|\nabla^2\E\LL(\upsilonv)-\nabla^2\E\LL(\upsilonvd)\|\le n\tilde L\|\upsilonv-\upsilonvd\|.
\end{EQA}

\subsubsection{Fisher theorem}
We controll the deviations of the maximizer of \(\LL\). The gradient reads
\begin{EQA}[c]
\nabla \LL(\upsilonv)=n\Xv-n\upsilonv+nf(\upsilonv)\frac{1}{2}\|\upsilonv\|\upsilonv+n\nabla f(\upsilonv)\|\upsilonv\|^{3}/3.
\end{EQA}
Setting this equal to zero we find that \(\tilde \upsilonv\) satisfies
\begin{EQA}[c]
\sqrt{n}\|\Xv-\tilde \upsilonv\|\le \frac{\sqrt n}{2}\|\tilde \upsilonv\|^2+\sqrt{n}\nabla f(\tilde\upsilonv)\|\tilde\upsilonv\|^{3}/3.
\end{EQA}
Using the fact that by Theorem 2.2 of \citep{AASP2013} \(\P(\tilde \upsilonv\in B_{2\sqrt{\frac{\dimn}{n}}}(0))\ge 1-2\ex^{-\dimtotal}\) such that \(\|\tilde \upsilonv\|\cong \sqrt{\dimn/n} \) and that \(\|\nabla f(\tilde\upsilonv)\|\cong \sqrt{n/\dimn}\) we obtain
\begin{EQA}[c]
\|\DPrc(\tilde \thetav-\thetavs)-\xivr\|^2=n\|\tilde \thetav-X_1\|^2\le n\|\tilde \upsilonv-\Xv\|^2\lesssim \dimn^2/n,
\end{EQA}
which shows that if \(\dimn^2/n\to 0\) we obtain the Fisher theorem.

\subsubsection{Wilks phenomenon}
Suppose for a moment that \(f\equiv 1\). One can see that the unique local maximizer \(\hat\upsilonv\) of
\begin{EQA}[c]
\hat\LL(\upsilonv)= n\mathbf X^\T \upsilonv-n\|\upsilonv\|^2/2+ n\|\upsilonv\|^{3}/3,
\end{EQA}
equals \(\lambda \Xv\) for some \(\lambda>0\) as only the term \(n\mathbf X^\T \upsilonv\) depends on the direction of \(\upsilonv\) and is maximized on balls with finite radius on the linear space spanned by \(\Xv\). We will show that \(\lambda=1+\delta(\Xv)\|\Xv\|\) where almost surely
\begin{EQA}[c]
\delta(\Xv)\to 4.
\end{EQA}
To see this note that the maximization problem reduces to solving
\begin{EQA}[c]
\argmax_{\lambda}\left\{ \lambda-\lambda^2/2+\|\Xv\|\lambda^3/3\right\}.
\end{EQA}
The solution can easily be obtained with first and second order criteria of maximality and is given as
\begin{EQA}
\lambda_{\max}&=&\frac{1-\sqrt{1-4\|\Xv\|}}{2\|\Xv\|}=\frac{4\|\Xv\|}{2\|\Xv\|(1+\sqrt{1-4\|\Xv\|})}\\
  &=&1+\frac{1-\sqrt{1-4\|\Xv\|}}{(1+\sqrt{1-4\|\Xv\|})}=1+\frac{4\|\Xv\|}{(1+\sqrt{1-4\|\Xv\|})^2}=:1+\tau(\Xv)\|\Xv\|.
\end{EQA}
Consequently \(\hat\upsilonv=(1+\tau(\Xv)\|\Xv\|)\Xv\).
Now if \(\hat\upsilonv\in \mathcal S\) this means that \(\tilde \upsilonv=\hat\upsilonv\) in our model, since for any other point \(\upsilonv\in\Upsilon\)
\begin{EQA}
\LL(\upsilonv)&=&n\mathbf X^\T \upsilonv-n\|\upsilonv\|^2/2+ f(\upsilonv)n\|\upsilonv\|^{3}\\
  &\le& n\mathbf X^\T \upsilonv-n\|\upsilonv\|^2/2+ n\|\upsilonv\|^{3}\\
  &\le& \max_{\upsilonv}\left\{n\mathbf X^\T \upsilonv-n\|\upsilonv\|^2/2+ n\|\upsilonv\|^{3}\right\}= \LL(\hat\upsilonv).
\end{EQA}
The event \(\{\hat\upsilonv\in \mathcal S\}\) is of strictly positive probability that depends on the choice of \(L>0\) and grows with \(n\to\infty\). Now observe that if \(\hat\upsilonv\in \mathcal S\)
\begin{EQA}
&&\nquad \Lr(\tilde \thetav)=\max_{\etav}\LL(\tilde\thetav,\etav)=\LL(\tilde \upsilonv)=\LL(\hat \upsilonv)\\
  &=&n\left((1+\tau(\Xv)\|\Xv\|)-(1+\tau(\Xv)\|\Xv\|)^2/2\right)\|\Xv\|^2+n(1+\tau(\Xv)\|\Xv\|)^3/3\|\Xv\|^3\\
  &=&n\|\Xv\|^2/2+n\|\Xv\|^3/3+n\left((\frac{1}{2}+\tau(\Xv))\|\Xv\|^4+\tau(\Xv)^2\|\Xv\|^5 +\tau(\Xv)^3\|\Xv\|^6/3\right).
\end{EQA}
By the definition of \(\Xv\) we have almost surely \(\lim n\|\Xv\|^2/\dimn \le \CONST\), such that if \(\dimn^2/n\to0\)
\begin{EQA}[c]
n\left((\frac{1}{2}+\tau(\Xv))\|\Xv\|^4+\tau(\Xv)^2\|\Xv\|^5 +\tau(\Xv)^3\|\Xv\|^6/3\right)=o_{\P}(1).
\end{EQA}
On the other hand we have due to \(f(0,\etav)=0\) for all \(\etav\in\R^{\dimn-1}\) with the orthogonal projection \(\Pi_{\etav}:\R^{\dimn}\mapsto \R^{\dimn-1}\) onto the nuisance component
\begin{EQA}[c]
\Lr(\thetavs)=\max_{\etav}\LL(\thetavs,\etav)=\max_{\etav}\left\{n\mathbf X^\T (0,\etav)-n\|\etav\|^2/2\right\}=n\|\Pi_{\etav}\mathbf X\|^2/2.
\end{EQA}
Consequently
\begin{EQA}
\Lr(\tilde \thetav)-\Lr(\thetavs)&=&n\|\Xv\|^2/2-n\|\Pi_{\etav}\mathbf X\|^2/2+n\|\Xv\|^{3}+o_{\P}(1)\\	
&=&nX_1^2/2+n\|\Xv\|^{3}+o_{\P}(1).
\end{EQA}
It is clear that if \(\dimn^3/n\rightarrow 0\) also \(n\|\Xv\|^{3}\rightarrow 0\) almost surely. Clearly \(nX_1^2\sim\chi^2_1\) for all \(n\in\N\). But if \(\dimn^3/n\nrightarrow 0\) obviously \(nX_1^2/2+n\|\Xv\|^{3}+o_{\P}(1)\) does not converge to a chisquare random variable with \(1\) degree of freedom. In consequence the Wilks phenomon does not occur on a set of positive probability if \(\dimn^3/n\nrightarrow 0\).

\bibliography{../../../sources/semiquellen}
  \end{document}